\documentclass[12pt,oneside,english]{amsart}
\usepackage[latin1]{inputenc}
\usepackage{amssymb}

\makeatletter
\newtheorem{theorem}{Theorem}
\newtheorem{lemma}{Lemma}

\newtheorem{definition}{Definition}
\newtheorem{condition}{Condition}
\newtheorem{proposition}{Proposition}
\newtheorem{remark}{Remark}
\newtheorem{conjecture}{Conjecture}

\newtheorem{corollary}{Corollary}

\numberwithin{equation}{section}
\usepackage{babel}

\makeatother
\begin{document}
\baselineskip=17pt

\title[Ramanujan and Labos primes]{ Ramanujan and Labos primes, their generalizations and classifications of primes}

\author{Vladimir Shevelev}
\address{Department of Mathematics \\Ben-Gurion University of the
 Negev\\Beer-Sheva 84105, Israel. e-mail:shevelev@bgu.ac.il}


\begin{abstract}
  Considering Ramanujan primes and the symmetric to them so-called Labos primes, we study their parallel properties , we study all primes with these properties (generalized Ramanujan and Labos primes) and construct two kinds of sieves for them. Finally, we give a further natural generalization of these constructions and pose some conjectures and open problems.\end{abstract}

\maketitle

\section{Introduction }

A very known Bertrand's postulate (1845) states that, for $x>1,$ always there exists a prime in interval $(x, 2x).$ This postulate very quickly-five years later- became a theorem due to Russian mathematician P.L.Chebyshev (cf., e.g., \cite{11}, Theorem 9.2).  In 1919, \enskip  S. Ramanujan \cite{9}-\cite{10} unexpectedly gave a new short and elegant proof of the Bertrand's postulate. In his proof appeared a sequence of primes
\begin{equation}\label{1.1}
2,11,17,29,41,47,59, 67, 71, 97, 101, 107, 127, 149, 151, 167,...
 \end{equation}
It is interesting that, for a long time this, important sequence was not presented in the Sloane's OEIS \cite{14}. Only in 2005\enskip  Sondow published it in OEIS (sequence A104272).
\begin{definition}\label{d1} For $n \geq 1,$ the \upshape nth Ramanujan  prime \slshape is the smallest positive integer  $(R_n)$  with the property that if $x\geq R_n,$  then $\pi(x) -\pi(x/2)\geq n.$
\end{definition}

In \cite{15}, Sondow obtained some estimates for $R_n$ and, in particular,  proved that, for every $n>1,\enskip R_n>p_{2n}.$ \enskip Laishram \cite{6} proved that $R_n<p_{3n}$ (a short proof of this result follows from a general Theorem 6 of the present paper, see Remark 2). Further, Sondow proved that, for $n\rightarrow\infty, \enskip R_n\sim p_{2n}.$ From this, denoting $\pi_R$ the counting function of the Ramanujan  primes, we have $R_{\pi_R(x)}\sim2\pi_R(x)\ln \pi_R(x).$
Since $R_{\pi_R(x)}\leq x<R_{\pi_{R(x)}+1},$
then $ x\sim p_{2\pi_R(x)}\sim2\pi_R(x)\ln \pi_R(x),$ as $x\rightarrow\infty,$ and we conclude that

 \begin{equation}\label{1.2}
\pi_R(x)\sim \frac {x} {2\ln x}\sim \frac {\pi(x)}{2}.
 \end{equation}
 Below we prove several other properties of the Ramanujan primes. An important role plays the following property.
\begin{theorem}\label{t1}
Let $p_n$ denote the $n$-th prime. If $p$ is an odd Ramanujan prime, such that $p_m< p/2< p_{m+1},$ then the interval  $(p,\enskip 2p_{m+1})$ contains a prime.
\end{theorem}
In 2003, Labos introduced the following sequence of primes (cf. \cite{14}, sequence A080359). We call them \slshape Labos primes \upshape, denoting $L_n$ the $n$-th Labos prime.
 \begin{definition}\label{d2} For $n \geq 1$, the \upshape nth Labos  prime\slshape \enskip is the smallest positive integer  $(L_n)$  for which $\pi(L_n)-\pi(L_n/2)=n.$
  \end{definition}
The first Labos primes are (see sequence A080359 in \cite{14}):
\begin{equation}\label{1.3}
2, 3, 13, 19, 31, 43, 53, 61, 71, 73, 101, 103, 109, 113, 139, 157, 173,...
 \end{equation}
Note that, since (\cite{14})
  \begin{equation}\label{1.4}
\pi(R_n)-\pi(R_n/2)=n,
 \end{equation}
 then, by the Definition 3, we have
\begin{equation}\label{1.5}
L_n\leq R_n.
 \end{equation}
 For them we prove a symmetric statement to Theorem 1.
\begin{theorem}\label{t2}
Let $p_n$ denote the $n$-th prime. If $p$ is an odd Labos prime, such that $p_m< p/2< p_{m+1},$ then the interval  $(2p_m,\enskip p)$ contains a prime.
\end{theorem}
It is clear that Theorems 1-2 are connected with some left-right symmetry in distribution of primes. Unfortunately, we cannot say about a \slshape precise \upshape left-right symmetry since the inequalities of type $R_1\leq L_1\leq R_2\leq L_2\leq...$ are broken from the very outset. Nevertheless, we prove the following theorem.
\begin{theorem}\label{t3}
If to consider \upshape all \slshape primes $\{R'_n\}$ and $\{L'_n\}$ for which Theorems 1-2 correspondingly are true, then for them we have
\begin{equation}\label{1.6}
R'_1\leq L'_1\leq R'_2\leq L'_2\leq...
 \end{equation}
\end{theorem}
Basing on this theorem, we give a natural simple classification of primes.\newline
\indent In this paper we prove Theorems 1-3 and several other properties of primes $\{R'_n\}$ and $\{L'_n\}$  and construct  two kinds of sieves for separation them among other primes.\newline
To obtain (\ref{1.6}), we should add to Ramanujan and Labos primes so-called pseudo-Ramanujan and pseudo-Labos primes.
 The first such primes are (see our sequence A164288 in \cite{14}):
\begin{equation}\label{1.7}
 109, 137, 191, 197, 283, 521, 617, 683, 907, 991, 1033, 1117, 1319,...
 \end{equation}
 \newpage
 The second ones are (A164294 in \cite{14}):
\begin{equation}\label{1.8}
 131, 151, 229, 233, 311, 571, 643, 727, 941, 1013, 1051, 1153, 1373,...
\end{equation}

\section{Proof of Theorems 1,2}

Below $p_n$ always denote the $n$-th prime. \newline
\indent We start with four conditions for odd primes.
\begin{condition}\label{co1}
Let $p=p_n,$ with  $n>1.$  Then all integers $(p+1)/2, (p+3)/2, ... , (p_{n+1}-1)/2 $ are composite numbers.
\end{condition}
\begin{condition}\label{co2}
Let, for an odd prime $p,$ we have $p_m< p/2< p_{m+1}.$ Then the interval  $(p,\enskip 2p_{m+1})$ contains a prime.
\end{condition} \begin{condition}\label{co3}
Let $p=p_n$ with $n\geq3.$ Then all integers $(p-1)/2, (p-3)/2, ... , (p_{n-1}+1)/2 $ are composite numbers.
\end{condition}
\begin{condition}\label{co4}
Let $p_m< p/2< p_{m+1}.$ Then the interval  $(2p_m,\enskip p)$ contains a prime.
\end{condition}
\begin{lemma}\label{l1}
 Conditions $\ref{co1}$ and $\ref{co2}$ are equivalent.
 \end{lemma}
 \bfseries Proof. \mdseries If Condition \ref{co1} is valid, then $p_{m+1}>(p_{n+1}-1)/2,$ i.e. $p_{m+1}\geq(p_{n+1}+1)/2.$ Thus $2p_{m+1}> p_{n+1}>p_n=p,$ and Condition \ref{co2} is valid; conversely, if Condition \ref{co2} satisfies, i.e. $p_{m+1}>p/2$
 and $2p_{m+1}>p_{n+1}>p=p_n.$  If $k$ is the least positive integer, such that $p_m<p_n/2<(p_n+k)/2<(p_{n+1}-1)/2$ and $(p_n+k)/2$ is prime, then $p_{m+1}=(p_n+k)/2$ and $p_{n+1}-1>p_n+k=2p_{m+1}>p_{n+1}.$ Contradiction shows that
 Condition \ref{co1} is valid.\newline
  $\blacksquare$\newline
  Quite analogously we obtain lemma on the equivalence of the second pair of the conditions.
\begin{lemma}\label{l2}
 Conditions $3$ and $4$ are equivalent.
 \end{lemma}
\indent Now we are able to prove Theorems 1,2.\newline
\indent In view of Lemma \ref{l1}, for proof of Theorem 1, it is sufficient to prove that, for Ramanujan primes, Condition \ref{co1} satisfies. If Condition \ref{co1} does not satisfy, then suppose that $p_m=R_n<p_{m+1}$ and  $k$ is the least positive integer, such that $q=(p_m+k)/2$ is prime not more than $(p_{m+1}-1)/2.$ Thus
  \begin{equation}\label{2.1}
  R_n=p_m<2q<p_{m+1}-1.
 \end{equation}
 From Definition \ref{d1} it follows that, $R_n-1$ is the maximal integer
 \newpage
 for which the equality
  \begin{equation}\label{2.2}
 \pi(R_n-1)-\pi((R_n-1)/2)=n-1
 \end{equation}
 holds. However, according to (\ref{2.1}), $\pi(2q)= \pi(R_n-1)+1$ and in view if the minimality of the prime $q,$ in the interval $((R_n-1)/2,q)$ there are not any prime. Thus $\pi(q)=\pi((R_n-1)/2)+1$ and
 $$ \pi(2q)-\pi(q)=\pi(R_n-1)-\pi((R_n-1)/2)=n-1.$$
 Since, by (\ref{2.1}), $2q>R_n,$ then this contradicts to the property of the maximality of $R_n$ in (\ref{2.2}). Thus the Theorem 1 follows.\newline
 \indent Theorem 2 is proved quite analogously, using Lemma 2.\newline$  \blacksquare$\newline
 \section{Pseudo-Ramanujan primes, over-Ramanujan primes and their Pseudo-Labos and over-Labos analogies}
 \begin{definition}\label{d3}
 Non-Ramanujan primes satisfying Condition $\ref{co1}\enskip ($or, equivalently, Condition $\ref{co2} )$ we call \upshape pseudo-Ramanujan primes, \slshape\enskip denoting the sequence of them  $\{R^*_n\}$ \upshape (see sequence $({1.7})).$
  \end{definition}

 \begin{definition}\label{d4}
 All primes satisfying Condition $\ref{co1}\enskip ($or, equivalently, Condition $\ref{co2} )$  we call \upshape over-Ramanujan primes, \slshape\enskip denoting the sequence of them $\{R'_n\}.$
  \end{definition}
  Note that all Ramanujan primes more then 2 are over-Ramanujan primes as well. Thus $R'_1=11.$
\newline
\indent Give a simple criterion for over-Ramanujan primes.
\begin{proposition}\label{p1}
$p_n\geq5$ is an over-Ramanujan prime if and only if $\pi(\frac {p_n} {2})=\pi(\frac {p_{n+1}} {2}).$
 \end{proposition}
   \bfseries Proof. \mdseries 1) Let $\pi(\frac {p_n} {2})=\pi(\frac {p_{n+1}} {2})$ is valid. From this it follows that if $p_k< p_n/2<p_{k+1},$ then between $p_n/2$ and $p_{n+1}/2$ there are not exist primes. Thus $p_{n+1}/2<p_{k+1}$ as well. Therefore, we have $2p_k< p_n<p_{n+1}<2p_{k+1},$ i.e. $p_n$ is an over-Ramanujan prime. Conversely, if $p_n$ is an over-Ramanujan prime, then $2p_k< p_n<p_{n+1}<2p_{k+1},$ and $\pi(\frac{p_n} {2})=\pi(\frac{p_{n+1}} {2})$ is valid. \newline
    $\blacksquare$
     \begin{definition}\label{d5}
 Non-Labos primes satisfying Condition $\ref{co1}\enskip ($or, equivalently, Condition $\ref{co2} )$ we call \upshape pseudo-Labos primes, \slshape\enskip denoting the sequence of them  $\{L^*_n\}$ \upshape (see sequence $({1.8})).$
  \end{definition}
  \newpage
  \begin{definition}\label{d6}
 All primes satisfying Condition $\ref{co3}\enskip ($or, equivalently, Condition $\ref{co4} )$  we call \upshape over-Labos primes, \slshape\enskip denoting the sequence of them $\{L'_n\}.$
  \end{definition}
   Note that all Labos primes more then 3 are over-Labos primes as well. Thus $L'_1=13.$
 Quite analogously to Proposition \ref{p1} we obtain the following criterion for over-Labos primes.
 \begin{proposition}\label{p2}
$p_n\geq5$ is over-Labos prime if and only if $\pi(\frac {p_{n-1}} {2})=\pi(\frac {p_{n}} {2}).$
 \end{proposition}
  \section{Proof of Theorem 3}
Now we prove a much stronger statement about the symmetry, which connected with the mutual behaviors of over-Ramanujan and over-Labos primes.\newline
\indent Note that the intervals of the form $(2p_m,\enskip 2p_{m+1})$ containing not more than one prime, contain neither over-Ramanujan nor over-Labos primes. Moving such intervals, consider the first from the remaining ones. The first its prime is  an over-Ramanujan prime $(R'_1).$ If it has only two primes, then the second prime is an over-Labos prime $(L'_1), $ and we see that $R'<L'; $ on the other hand if it has $k$ primes, then beginning with the second one and up to the $(k-1)$-th we have primes which are simultaneously over-Ramanujan and over-Labos primes. Thus, taking into account that the last prime is only over-Labos prime , we have
 $$R'_1<L'_1=R'_2=L'_2=R'_3=...=L'_{k-1}=R'_{k-1}<L'_{k}.$$
 The second remaining interval begins with an over-Ramanujan prime and the process repeats. \newline $\blacksquare$
  \section{Prime gaps}
 Recently, Sondow, Nicholson and Noe \cite{16} showed that, if to consider a run of consecutive Ramanujan primes $p=R_l,...,q=R_k,$ then the interval $\frac{1}{2}(p+1),\enskip \frac{1}{2}(q+1)$ is free from primes. According their definition, the interval $[a,b]$ is a prime gap, if none of the numbers $a,a+1,...,b$ is prime. Nevertheless, their result is far from a complete characterization of the prime gaps in a usual sense. E.g., we have a run $\{2521,2531\}$ of consecutive Ramanujan primes which gives a "prime gap" $[\frac{2521+1}{2},\enskip \frac{2531+1}{2}]=[1261,1266].$ However, the real prime gap is much more: (1259,1277). A better result one can obtain using over-Ramanujan primes. Indeed, the used in \cite{16} properties of Ramanujan primes are valid for all over-Ramanujan primes, while runs of consecutive over-Ramanujan primes, generally speaking, are longer. 
 \newpage 
 E.g., instead of run $\{2521,2531\}$ of Ramanujan primes, we have run $\{2521,2531,2539,2543,2549\}$ of over-Ramanujan primes. This gives the interval $[\frac{2521+1}{2},\enskip \frac{2549+1}{2}]=[1261,1275]$ which is free from primes and very close to the real gap. In general, since  over-Ramanujan primes satisfy Condition 1, then to every run of consecutive over-Ramanujan primes $p=R'_l,...,q=R'_k$ corresponds interval $[\frac{1}{2}(p+1),\enskip \frac{1}{2}(q+1)]$ which contains no primes. Note that the next after $q$ prime $q'$ gives an additional improvement of lower estimate of size $(L)$ of the considered prime gap. Indeed, we know that $q'$ is necessarily an over-Labos prime. Since the over-Labos primes satisfy Condition 3, then all numbers $\frac{q'-1}{2}, \enskip \frac{q'-2}{2},...,\frac {q+1}{2}$ are composite. Hence $L\geq\frac{q'-p}{2}.$
For example, consider the run $\{227,229,233,239,241\}$ of over-Ramanujan primes (all these primes are Ramanujan). The following prime is $q'=251.$ Thus, for the gap containing $(227+1)/2=114,$ we have $L\geq\frac {251-227}{2}=12$ (the exact value of $L$ here is 14).
   \section{The first sieve for selection of the over-Ramanujan primes from all primes}

\indent  Recall that the Bertrand sequence $\{b(n)\}$ is defined as $ b(1)=2,$ and, for $n\geq2,\enskip b(n)$ is the largest prime less than $2b(n-1)$ (see A006992 in \cite{14}):
\begin{equation}\label{6.1}
2, 3, 5, 7, 13, 23, 43, 83, 163, 317, 631, 1259, 2503, 5003,...
 \end{equation}
Put
\begin{equation}\label{6.2}
B_0=\{b^{(0)}(n)\}=\{{b(n)}\}.
 \end{equation}
Further we build sequences $B_1=\{b^{(1)}(n)\}, B_2=\{b^{(2)}(n)\},...$ according the following inductive rule: if we have sequences $B_0,...,B_{k-1},$ let us consider the minimal prime $p^{(k)}\not \in \bigcup_{i=1}^{k-1}B_i.$ Then the sequence $\{b^{(k)}(n)\}$ is defined as $b^{(k)}(1)=p^{(k)},$ and, for $n\geq2,\enskip b^{(k)}(n) $ is the largest prime less than $2b^{(k)}(n-1).$ So, we obtain consequently:
\begin{equation}\label{6.3}
B_1=\{11, 19, 37, 73,...\}
 \end{equation}
 \begin{equation}\label{6.4}
B_2=\{17, 31, 61, 113,...\}
 \end{equation}
 \begin{equation}\label{6.5}
B_3=\{29, 53, 103, 199,...\}
 \end{equation}
 etc., such that, putting $p^{(1)}=11,$ we obtain the sequence
 \begin{equation}\label{6.6}
\{p^{(k)}\}_{k\geq1}=\{11, 17, 29, 41, 47, 59, 67, 71, 97, 101, 107, 109, 127, ...\}
 \end{equation}
 Sequence (\ref{6.6}) coincides with sequence (\ref{1.1}) of the Ramanujan primes from the second and up to the $12$-th term, but the $13$-th term of this sequence is $109 $ which is the first term of sequence (\ref{1.7}) of the pseudo-Ramanujan 
 \newpage
 primes.
 
 \begin{theorem}\label{t4} For $n\geq1,$ we have
 \begin{equation}\label{6.7}
p^{(n)}=R'_n.
 \end{equation}
 \end{theorem}
\bfseries Proof. \mdseries The least omitted prime in (\ref{6.1}) is $p^{(1)}=11=R'_1$;  the least omitted prime in the union of (\ref{6.2}) and (\ref{6.3}) is $p^{(2)}=17=R'_2.$ We use the

induction. Let we have already built primes $$p^{(1)}=11, p^{(2)},...,p^{(n-1)}=R'_{n-1}.$$
 Let $q$ be the least prime which is omitted in the union $\bigcup_{i=1}^{n-1}B_i,$  such that $q/2$ is in interval $(p_{m}, p_{m+1}).$ According to our algorithm, $q$ which is dropped should not be the largest  prime in the interval $(p_{m+1}, 2p_{m+1}).$ Then there are primes in the interval  $q, 2p_{m+1})$; let $r$ be one of them. Then we have $2p_m<q<r<2p_{m+1}.$ This means that $q,$ in view of its minimality between the dropping primes which are more than $R'_{n-1}=p^{(n-1)},$ is the least over-Ramanujan prime more than $R'_{n-1}$ and  the least prime of the form $p^{(k)}$ more than $p^{(n-1)}$. Therefore, $q=p^{(n)}=R'_{n}.\newline\blacksquare$\newline
 \indent Quite analogously, using sequence $ c(1)=2,$ such that for $n\geq2,\enskip c(n)$ is the smallest prime more than $2c(n-1)$ (see A055496 in \cite{14}), one can construct a sieve for over-Labos primes.

 \section{The second sieve for selection of the over-Ramanujan primes from all primes}
 The theorem on precise symmetry allows construct the second  sieve for over-Ramanujan primes. \newline
 \indent Consider consecutive  intervals of the form $(2p_n,\enskip 2p_{n+1}),\enskip n=1,2,...$ Remove all of them which contain less than two primes. For every remain interval, we write primes (in increasing order) except of the last one. Then all remain primes are over-Ramanujan.\newline
 Let us demonstrate this sieve. For primes 2,3,5,7,11,... consider intervals
 \begin{equation}\label{7.1}
  (4,6),(6,10),(10,14),(14,22),(22,26),(26,34), (34,38),(38,46),...
  \end{equation}
Remove those of them which contain less than two primes. We have the following sequence of intervals:
\begin{equation}\label{7.2}
 (10,14),(14,22),(26,34,),(38,46),(46,58),(58,62),(62,74),...
  \end{equation}
 Now we write all primes from these intervals, excluding the \slshape last \upshape \enskip primes. Then we obtain sequence (\ref{6.6}).
  \newpage
   \indent Quite analogously we obtain the second sieve for over-Labos primes. This sequence we can obtain in a parallel way. It is sufficient to write all primes from the last sequence of intervals, excluding the \slshape first \upshape \enskip primes. Thus we obtain the sequence (cf. A164333 in \cite{14})
 \begin{equation}\label{6.3}
13, 19, 31, 43, 53, 61, 71, 73, 101, 103, 109, 113, 131, 139, 151, 157,...
 \end{equation}

\section{A classification of primes}

\indent In connection with the considered construction, let us consider the following classification of primes.\newline \indent 1) Two first primes 2,3 form a separate set of primes. \newline
\indent 2) If $p\geq11$ is an over-Ramanujan but over-Labos prime, then, in connection with the second sieve, we call $p$ a \slshape left prime\upshape \enskip (cf. A166307 in \cite{14}):
$$11, 17, 29, 41, 47, 59, 67, 97, 107, 127, 137, 149, 167, 179, 197, 227,...  $$
\indent 3) If $p\geq5$ is an over-Labos but over-Ramanujan prime, then we call $p$ a \slshape right prime.\upshape. \newline The first terms of this sequence are
$$13, 19, 31, 43, 53, 61, 73, 103, 113, 131, 139, 157, 173, 193, 199, 251,...  $$
\indent 4) If $p$ is simultaneously over-Ramanujan and over-Labos prime, then we  call it \slshape a central prime \upshape \enskip ( sequence A166252 in \cite{14}):
$$71, 101, 109, 151, 181, 191, 229, 233, 239, 241, 269, 283, 311, 349,...  $$
\indent 5) Finally, the rest primes it is natural to call \slshape isolated primes \upshape \enskip (sequence A166251 in \cite{14}):
$$5, 7, 23, 37, 79, 83, 89, 163, 211, 223, 257, 277, 317, 331, 337, 359,...  $$
 Note that from the second sieve the following result follows.
 \begin{proposition}\label{p3}
Let $l_n,\enskip r_n$ denote the $n$-th left prime and the $n$-th right prime correspondingly. Then, for $n\rightarrow\infty,$ we have $l_n\sim r_n.$
 \end{proposition}
 \bfseries Proof.\mdseries \enskip Beginning with Hoheisel's proof \cite{5} that, for $x>x_0(\varepsilon),$ the interval $(x, x+x^{1-\frac {1} {33000}+\varepsilon}]$ always contains a prime (with numerous improvements up to currently the best result of Baker, Harman and Pintz [2]), it is known that $p_{n+1}-p_n=o(p_n).$ Since, by the construction, $l_n$ and $r_n$ belong to the same interval of the form $(2p_{m(n)},\enskip 2p_{m(n)+1}),$ then $r_n-l_n<2(p_{m(n)+1}-p_{m(n)})=o(l_n),$ then the statement follows.\newline$\blacksquare$
\newpage

\section{On density of over-Ramanujan and over-Labos primes}
\indent Unfortunately, the research of the obtained two kinds of sieves seems much more difficult than the research of the Eratosthenes one for primes. Nevertheless, some very simple probabilistic arguments lead to a very plausible conjecture about the density of over-Ramanujan and over-Labos primes.
\begin{conjecture}\label{con1}
Let $\pi_{R'}(x)$ be the counting function of over-Ramanujan numbers not exceeding $x.$ Then
 \begin{equation}\label{9.1}
\pi_{R'}(x)\sim (\frac{1}{2}+\frac{1}{e^2-1})\pi(x)=0.6565176...\pi(x).
 \end{equation}
\end{conjecture}
\bfseries Heuristic proof.\mdseries \enskip Consider asymptotically $\frac{\pi(n)}{2}$ intervals of the form $(2p_m,\enskip 2p_{m+1})$ covering all $\pi(n)$ primes. Berend \cite{3} noticed that the number of primes which are not $v$-over-Ramanujan among $\pi(n)$ primes exactly equals to the number of the considered intervals containing at least one prime. Indeed, a prime is not over-Ramanujan if and only if it is the last prime of a such interval. It is well known (\cite{8}) that for large $n$ a random interval between two consecutive primes has length $\ln p_n.$ Thus a random interval of the considered form has length $2\ln p_n$ and, according to the Cram\'{e}r model \cite{4}, the number of primes in a such random interval has the binomial $(2\ln p_n,\enskip \frac {1} {\ln p_n})$ distribution which, for large $\ln p_n,$ has a good approximation by a Poisson distribution with parameter $\lambda=2.$ Let us calculate asymptotically the number of over-Ramanujan primes not exceeding $n,$ using this model. A random interval contains $k$ primes with the probability $\mathrm{P}[X=k]=\frac {2^k}{k!}e^{-2},\enskip k=0,1,2...\enskip .$ On the other hand, an interval contains $i\geq1$ over-Ramanujan primes if and only if it contains $i+1$ primes. It is clear that we consider a random interval in the condition that it already contains a prime. Thus the total number of
of over-Ramanujan primes not exceeding $n$ asymptotically equals to
$$\frac{\pi(n)}{2}\sum_{i\geq1}\mathrm{P}[X=i+1 | X\geq1]=$$ $$\frac{\pi(n)}{2}(1-e^{-2})^{-1}\sum_{i\geq1}\frac{2^{i+1}}{(i+1)!}e^{-2}i=$$ $$\pi(n)\frac {e^{-2}}{1-e^{-2}}\sum_{i\geq1}(\frac{2^{i}}{i!}-\frac{2^{i}}{(i+1)!}  )=$$$$\pi(n)\frac{e^{-2}}{1-e^{-2}}((e^2-1)-\frac{1}{2}(e^2-3))=  $$ $$\frac{\pi(n)}{2}\frac{e^{-2}}{1-e^{-2}}(e^2+1) $$
\newpage
and ({9.1}) follows. $\enskip \blacksquare$\newline
\indent Greg Martin \cite{7} did the corresponding calculations for the first million primes $p$, and found that for approximately $61.2\% $ of them have a prime in the interval $(p, 2p_{n+1}).$ Since in this case $\ln p_n$ is small (less than 17), then an error of order $4\% $ is quite acceptable. Moreover, Martin conjectured that the probability is $\frac{2}{3}.$ This differs from the probability ({9.1}) only on $1\%!$\newline
\indent Note that, if Conjecture \ref{con1} is true, then, using (\ref{1.2}), for the counting function $\pi_{R^*}(x)$ of pseudo-Ramanujan primes, we have
 \begin{equation}\label{9.2}
\pi_{R^*}(x)\sim \frac {\pi(x)}{e^2-1}\pi(x)=0.15651...\pi(x),
 \end{equation}
 such that the proportion of Ramanujan primes among all over-Ramanujan primes is approximately 0.76159...\newline
 Using Theorem \ref{t3}, note that, if Conjecture \ref{con1} is true, then, for the counting function $\pi_{L'}(x)$ of over-Labos primes, we have
  \begin{equation}\label{9.3}
\pi_{L'}(x)\sim \pi_{R'}(x)\sim (\frac{1}{2}+\frac{1}{e^2-1})\pi(x).
 \end{equation}
 Show that events $R':$ "a prime is over-Ramanujan" and $L':$ "a prime is over-Labos" are independent. Indeed, denoting events $r:$ "a prime is right", $l:$ "a prime is left" and $Is:$ "a prime is isolated", we have
 \begin{equation}\label{9.4}
 \mathrm{P}[R' | L']=1-\mathrm{P}[l]-\mathrm{P}[Is]; \enskip \mathrm{P}[R' |\overline{L'} ]=1-\mathrm{P}[r]-\mathrm{P}[Is].
 \end{equation}
 Hence, in view of $\mathrm{P}[l]=\mathrm{P}[r]$ (cf. Proposition 3), we have
 \begin{equation}\label{9.5}
 \mathrm{P}[R' | L']=\mathrm{P}[R' |\overline{L'} ].
 \end{equation}
 \indent Therefore, if Conjecture \ref{con1} is true, then, for the counting function $\pi_l(x),\enskip \pi_r(x),\enskip \pi_c(x)$ and $\pi_{is}(x) $ of the left, right, central and isolated primes correspondingly of our classification of primes, we have
  \begin{equation}\label{9.6}
\pi_l(x)\sim\pi_r(x)\sim (\frac{1}{4}-\frac{1}{(e^2-1)^2})\pi(x)=0.2255...\pi(x),
 \end{equation}
  \begin{equation}\label{9.7}
\pi_c(x)\sim(\frac{1}{2}+\frac{1}{e^2-1})^2\pi(x)=0.4310...\pi(x),
 \end{equation}
 \begin{equation}\label{9.8}
  \pi_{is}(x)\sim (\frac{1}{2}-\frac{1}{e^2-1})^2=0.1179...\pi(x),
 \end{equation}
 such that $\pi_r(x)+\pi_l(x)+\pi_c(x)+\pi_{is}(x)=\pi(x).$
\section{A generalization}

Let us consider a natural generalization of Ramanujan primes.
\begin{definition}\label{d7} For a given real $v>1,$ we call a \;\upshape $v$-Ramanujan prime $R_v(n),\enskip n\geq1,$\slshape\enskip the smallest integer with the property: if  $x \geq R_v(n),$ then
\newpage 
 $\pi(x) - \pi(x/v) \geq n.$ \end{definition}
\indent It is easy to see that $R_v(n)$ is indeed a prime. Moreover, as in \cite{15}, one can prove that
\begin{equation}\label{10.1}
R_v(n)\sim p_{((v/(v-1))n)},
\end{equation}
 as n tends to the infinity. Let $\pi_R^{(v)}(x)$ be the counting function of $v$-Ramanujan primes. Then we have (cf. (\ref{1.2}))
 \begin{equation}\label{10.2}
 \pi_R^{(v)}(x)\sim(1-1/v)\pi(x).
\end{equation}
Put
 \begin{equation}\label{10.3}
 \kappa(v)=\begin{cases}0,& if \;\; v\enskip is\enskip not \enskip ratio \enskip of \enskip primes;
\\r,& if\;\;v=\frac{r}{q},\enskip where\enskip r\enskip and\enskip q\enskip are\enskip primes.\end{cases}
\end{equation}
The following theorem is proved by the same way as Theorem 1.
\begin{theorem}\label{t5}
Let $p_n$ denote the $n$-th prime. Let $v>1$ be a given real number. If $p>\max(2v,\enskip \kappa(v))$ is an $v$-Ramanujan prime, such that $p_m< p/v< p_{m+1},$ then the interval  $(p,\enskip \lceil vp_{m+1}\rceil+\varepsilon)$ contains a prime.
\end{theorem}
\begin{remark}
Condition $p>\max(2v,\enskip \kappa(v))$ allows to avoid the cases $p=2v$ and $p=vq$ with a prime $q,$ when the condition $p_n< p/v< p_{n+1}$ is impossible.
\end{remark}
Let us find an upper estimate for the $n$-th $v$-Ramanujan prime.
\begin{theorem}\label{t6} If $n\geq\frac{1}{k}\max (6k,\enskip e^v,\enskip v^{(0.79677\frac{k-1}{k}v-1)^{-1}}),$ then, for $v\geq1.25507\frac{k}{k-1},$ we have
\begin{equation}\label{10.4}
R_v(n)\leq p_{kn}.
\end{equation}
\end{theorem}
\bfseries Proof.\mdseries\enskip It is sufficient to show that $\pi(\frac{p_{kn}}{v})\leq(k-1)n.$ Indeed, then we have
$\pi(p_{kn})-\pi(\frac{p_{kn}}{v})\geq kn-(k-1)n=n.$ We use the following known results ([1], \cite {12}-\cite{13}):
\begin{equation}\label{10.5}
p_n<n\ln n+n\ln\ln n,\enskip n\geq6;
\end{equation}
\begin{equation}\label{10.6}
p_n>n\ln n;
\end{equation}
\begin{equation}\label{10.7}
\pi(x)<1.25506 \frac {x} {\ln x},\enskip x>1.
\end{equation}
Note that, $\frac {p_{kn}}{v}>\frac {kn}{v}>\frac{kn}{e^v}.$  Hence, by the condition, $\frac{p_{kn}}{v}>1.$
 By (\ref{10.5})-(\ref{10.7}), we have
$$\pi(\frac{p_{kn}}{v})<1.25506\frac {p_{kn}}{v \ln (\frac {p_{kn}}{v})}<$$
 \newpage
 $$1.25506\frac{kn}{v}\cdot\frac {\ln(kn)+\ln(\ln(kn))}{\ln (\frac {kn\ln(kn)}{v})}=$$
 $$1.25506\frac{kn}{v}(1+\frac {\ln v}{\ln (\frac {kn\ln(kn)}{v})}).$$
Taking into account that, by the condition, $\ln (kn)>v,$ we have
$$\pi(\frac{p_{kn}}{v})<1.25506\frac{kn}{v}(1+\frac {\ln v}{\ln (kn)}).$$
Finally, note that, by the condition, $\frac {\ln v}{\ln (kn)}\leq 0.7968\frac{k-1}{k}v-1.$ Therefore,
$$\pi(\frac{p_{kn}}{v})<1.25506\cdot0.79677(k-1)n<(k-1)n.\enskip \blacksquare$$
\begin{corollary}
\begin{equation}\label{10.8}
R_3(n)<p_{2n},\enskip n\geq1.
\end{equation}
\begin{equation}\label{10.9}
R_{1.8}(n)<p_{4n},\enskip n\geq1.
\end{equation}
\end{corollary}
\bfseries Proof.\mdseries \enskip By Theorem \ref{t6}, for $v=3,\enskip k=2,$  we get the required inequality for $n\geq279.$ Using a small computer verification for $n<279,$ we obtain ({10.8}). In case  $v=1.8,\enskip k=4,$ we get the required inequality for $n\geq2370.$ Using a simple computer verification for $n<2370,$ we obtain ({10.9}).$\enskip \blacksquare$\newline
\indent The first terms of $1.8$-Ramanujan primes are
\begin{equation}\label{10.10}
2,11,17,37,43,59,61,79,97,101,103,137,163,167,191,211,...
\end{equation}
\begin{remark}\label{r2}
In case $v=2,\enskip k=3,$ by Theorem $6,$ we find $R_n=R_2(n)<p_{3n},$ for $n\geq22398.$ A simple computer verification for $n<22398$ leads to the Laishram result \cite{6}.
\end{remark}
\begin{definition}\label{d8}
 Every prime $p>\max(2v,\enskip \kappa(v))$ is called a \upshape $v$-over-Ramanujan prime,\slshape\enskip if, as soon as $p_m< p/v< p_{m+1},$ then the interval  $(p,\enskip vp_{m+1})$ contains a prime.
 \end{definition}
\begin{definition}\label{d9}
 A $v$-over-Ramanujan but $v$-Ramanujan prime is called \upshape $v$-pseudo-Ramanujan prime.
 \end{definition}
\slshape Now $v$-Labos primes, $v$-over-Labos primes and $v$-pseudo-Labos prime are introduced quite symmetrically \upshape \enskip(cf. Section 3). In particular, the following statements are valid.
\begin{theorem}\label{t7}
Let $p_n$ denote the $n$-th prime. Let $v>1$ be a given real number. If $p>\max(2v,\enskip \kappa(v))$ is an $v$-Labos prime, such that $p_m< p/v< p_{m+1},$ then the interval  $(\lfloor vp_m\rfloor-\varepsilon,\enskip p)$ contains a prime.
\end{theorem}
\newpage
\begin{theorem}\label{t8}
For sequences $\{R'_v(n)\}$ and $\{L'_v(n)\}$ of $v$-over-Ramanujan and $v$-over-Labos primes, we have
\begin{equation}\label{10.11}
R'_v(1)\leq L'_v(1)\leq R'_v(2)\leq L'_v(2)\leq...
 \end{equation}
\end{theorem}
A generalization of the first sieve for $v$-over-Ramanujan primes, $v\geq2,$ is based on Bertrand-like sequence $\{b_v(n)\}$ which is defined as $ b_v(1)=2,$ and, for $n\geq2,\enskip b_v(n)$ is the largest prime less than $\lceil vb_v(n-1)\rceil+\varepsilon.$ A generalization of the second sieve for $v$-over-Ramanujan primes is based on the sequence of intervals
\begin{equation}\label{10.12}
(\lfloor2v\rfloor-\varepsilon,\lceil3v\rceil+\varepsilon),(\lfloor3v\rfloor-\varepsilon,\lceil5v\rceil+\varepsilon),
(\lfloor5v\rfloor-\varepsilon,\lceil7v\rceil+\varepsilon),...
 \end{equation}
 with the removing the intervals containing less than two primes (cf. (\ref{7.2})).  For every remain interval, we write primes (in increasing order) except of the last one. Then all remain primes are $v$-over-Ramanujan. \newline
 \indent For example, if $v=3,$ we obtain the following sequence of $3$-over-Ramanujan primes (cf. sequence A164952 in \cite{14}):
 \begin{equation}\label{10.13}
2, 3, 11, 17, 23, 29, 41, 43, 59, 61, 71, 73, 79, 97, 101, 103, 107,...
 \end{equation}
 Furthermore, one can obtain a $v$-classification of primes, including $v$-left, $v$-right, $v$-central and $v$-isolated primes (cf. Section 8). In particular, if $l_v(n),\enskip r_v(n)$ denote the $n$-th $v$-left prime and the $n$-th $v$-right prime correspondingly, then, for $n\rightarrow\infty,$ we have $l_n\sim r_n.$\newline
\indent A generalization of Conjecture \ref{con1} (with the quite similar heuristic proof) is the following.
\begin{conjecture}\label{con2}
Let $\pi_{R_v'}(x)$ be the counting function of $v$-over-Ramanujan numbers not exceeding $x.$ Then
 \begin{equation}\label{10.14}
\pi_{R_v'}(x)\sim (\frac{v-1}{v}+\frac{1}{e^{v}-1})\pi(x).
 \end{equation}
\end{conjecture}
\indent Note that, if Conjecture \ref{con2} is true, then, using (\ref{10.2}), for the counting function $\pi_{R_v^*}(x)$ of $v$-pseudo-Ramanujan primes, we have
 \begin{equation}\label{10.15}
\pi_{R_v^*}(x)\sim \frac {\pi(x)}{e^{v}-1},
 \end{equation}
 such that the proportion of $v$-pseudo-Ramanujan primes among all $v$-over-Ramanujan primes is
 $\frac{v}{(v-1)e^v+1}.$ Thus this proportion tends to 1, when $v\rightarrow1,$

  and fast tends to 0, when $v\rightarrow\infty.$\newline

 Using Theorem \ref{t8}, note that, if Conjecture \ref{con2} is true, then, for the counting function $\pi_{L_v'}(x)$ of over-Labos primes, we have
  \begin{equation}\label{10.16}
\pi_{L_v'}(x)\sim \pi_{R_v'}(x)\sim (\frac{v-1}{v}+\frac{1}{e^{v}-1})\pi(x).
 \end{equation}
 \newpage
 \indent Furthermore, if Conjecture \ref{con2} is true, then, for the counting function $\pi_{l_v}(x),\enskip \pi_{r_v}(x),\enskip \pi_{c_v}(x)$ and $\pi_{is_v}(x) $ of the $v$-left, $v$-right, $v$-central and $v$-isolated primes correspondingly of our $v$-classification of primes, we have
  \begin{equation}\label{10.17}
\pi_{l_v}(x)\sim\pi_{r_v}(x)\sim (\frac{v-1}{v}+\frac{1}{e^{v}-1})(\frac{1}{v}-\frac{1}{e^{v}-1})\pi(x),
 \end{equation}
  \begin{equation}\label{10.18}
\pi_{c_v}(x)\sim (\frac{v-1}{v}+\frac{1}{e^{v}-1})^2\pi(x),
 \end{equation}
 \begin{equation}\label{10.19}
   \pi_{is_v}(x)\sim (\frac{1}{v}-\frac{1}{e^{v}-1})^2\pi(x),
 \end{equation}
such that $\pi_{r_v}(x)+\pi_{l_v}(x)+\pi_{c_v}(x)+\pi_{is_v}(x)=\pi(x).$ \newline
\indent Note that, by the heuristic arguments, the approximations given by formulas (\ref{10.14})-(\ref{10.19}) should be considered only for large magnitudes of $v\ln p_n.$

\section{Other open problems}

1) For $v>1$ to estimate the smallest pseudo-$v$-Ramanujan prime.\newline
\indent 2) For $v>1$ to estimate the smallest $v$-central prime.\newline
\indent 3) For $v>1$ to estimate the smallest $v$-isolated prime.

\section{Acknowledgments}
The author is grateful to Daniel Berend (Ben Gurion University, Israel) and Greg Martin (University of British Columbia, Canada) for important private communications \cite{3}, \cite{7} and very useful discussions.

\end{document}